\documentclass[conference]{HPCSTrans}
\IEEEoverridecommandlockouts
\usepackage{amsmath,amssymb,enumerate}
\usepackage{cite,graphics,graphicx}
\usepackage{amsmath,amssymb,amsthm,algorithm}
\usepackage{listings,float}
\usepackage{hyperref}
\usepackage{multirow}
\usepackage{pbox}
\usepackage{rotating}
\usepackage{array}
\usepackage[T1]{fontenc}

\begin{document}
%
% paper title
% can use linebreaks \\ within to get better formatting as desired
\title{
\vspace{-1em}
Exponential Integrators on Graphic Processing Units
}

% author names and affiliations
% use a multiple column layout for up to three different
% affiliations
\author{\IEEEauthorblockN{Lukas Einkemmer, Alexander Ostermann}
\IEEEauthorblockA{
\textit{Department of Mathematics}\\
\textit{University of Innsbruck}\\
\textit{Innsbruck, Austria}\\
\textit{lukas.einkemmer@uibk.ac.at, alexander.ostermann@uibk.ac.at} \\
\\
{\small \copyright 2013 IEEE. Personal use of this material is permitted. Permission from IEEE must be obtained for all other uses, in any current or} \\
{\small future media, including reprinting/republishing this material for advertising or promotional purposes, creating new collective works, } \\
{\small for resale or redistribution to servers or lists, or reuse of any copyrighted component of this work in other works.}
% http://www.ieee.org/publications_standards/publications/rights/rights_policies.html
%{\small
%\copyright 2013 IEEE. Personal use of this material is permitted. Permission from IEEE must be obtained for all other uses, in any current or future media, including reprinting/republishing this material for advertising or promotional purposes, creating new collective works, for resale or redistribution to servers or lists, or reuse of any copyrighted component of this work in other works.
%l}
\vspace{-2em}
}~\thanks{}

%\and
%\IEEEauthorblockN{Alexander Ostermann}
%\IEEEauthorblockA{
%\textit{Department of Mathematics}\\
%\textit{University of Innsbruck}\\
%\textit{Innsbruck, Austria}\\
%\textit{alexander.ostermann@uibk.ac.at}
%}
}

% conference papers do not typically use \thanks and this command
% is locked out in conference mode. If really needed, such as for
% the acknowledgment of grants, issue a \IEEEoverridecommandlockouts
% after \documentclass

% for over three affiliations, or if they all won't fit within the width
% of the page, use this alternative format:
%
%\author{\IEEEauthorblockN{Michael Shell\IEEEauthorrefmark{1},
%Homer Simpson\IEEEauthorrefmark{2},
%James Kirk\IEEEauthorrefmark{3},
%Montgomery Scott\IEEEauthorrefmark{3} and
%Eldon Tyrell\IEEEauthorrefmark{4}}
%\IEEEauthorblockA{\IEEEauthorrefmark{1}School of Electrical and Computer Engineering\\
%Georgia Institute of Technology,
%Atlanta, Georgia 30332--0250\\ Email: see http://www.michaelshell.org/contact.html}
%\IEEEauthorblockA{\IEEEauthorrefmark{2}Twentieth Century Fox, Springfield, USA\\
%Email: homer@thesimpsons.com}
%\IEEEauthorblockA{\IEEEauthorrefmark{3}Starfleet Academy, San Francisco, California 96678-2391\\
%Telephone: (800) 555--1212, Fax: (888) 555--1212}
%\IEEEauthorblockA{\IEEEauthorrefmark{4}Tyrell Inc., 123 Replicant Street, Los Angeles, California 90210--4321}}

% use for special paper notices
%\IEEEspecialpapernotice{(Invited Paper)}

% no paragraph intentation
%\setlength{\parindent}{0in}

% make the title area

\maketitle

\begin{abstract}
%\boldmath
In this paper we revisit stencil methods on GPUs in the context of exponential integrators. We further discuss boundary conditions, in the same context, and show that simple boundary conditions (for example, homogeneous Dirichlet or homogeneous Neumann boundary conditions) do not affect the performance if implemented directly into the CUDA kernel. In addition, we show that stencil methods with position-dependent coefficients can be implemented efficiently as well.
                As an application, we discuss the implementation of exponential integrators for different classes of problems in a single and multi GPU setup (up to 4 GPUs). We further show that for stencil based methods such parallelization can be done very efficiently, while for some unstructured matrices the parallelization to multiple GPUs is severely limited by the throughput of the PCIe bus.
\end{abstract}
% IEEEtran.cls defaults to using nonbold math in the Abstract.
% This preserves the distinction between vectors and scalars. However,
% if the conference you are submitting to favors bold math in the abstract,
% then you can use LaTeX's standard command \boldmath at the very start
% of the abstract to achieve this. Many IEEE journals/conferences frown on
% math in the abstract anyway.
% no keywords
\vspace{0.1in}

% For peer review papers, you can put extra information on the cover
% page as needed:
% \ifCLASSOPTIONpeerreview
% \begin{center} \bfseries EDICS Category: 3-BBND \end{center}
% \fi
%
% For peerreview papers, this IEEEtran command inserts a page break and
% creates the second title. It will be ignored for other modes.
\IEEEpeerreviewmaketitle

\section{INTRODUCTION}

The emergence of graphic processing units as a massively parallel computing architecture as well as their inclusion in high performance computing systems have made them an attractive platform for the parallelization of well established computer codes.

Many problems that arise in science and engineering can be modeled as differential equations. In most circumstances the resulting equations are sufficiently complex such that they can not be solved exactly. However, an approximation computed by the means of a given numerical scheme can still serve as a valuable tool for scientists and engineers. The collection of techniques generally referred to as general-purpose computing on graphics processing units (GPGPU) provide the means to speed up such computations significantly (see e.g. \cite{murray2012} or \cite{Micikevicius2009}).

%In order to simulate a continuous problem on a computer a discretization procedure has to be employed. The method of lines, for example, transforms a partial differential equation into a problem that is discretized in space but still has a continuous dependence on time. The resulting ordinary differential equation is then solved by a time integration scheme.
If a finite difference approximation in space is employed (such methods are widely used in computational fluid dynamics, for example), stencil methods provide an alternative to storing the matrix in memory (see e.g. \cite{Berkeley2009}). In many instances, this is advantageous both from a memory consumption as well as from a performance standpoint. The resulting system of ordinary differential equations then has to be integrated in time.

Much research has been devoted to the construction of efficient time integration schemes as well as their implementation (see e.g. \cite{hairer1} and \cite{hairer2}). The implementation of Runge--Kutta methods, which are the most widely known time integration schemes, on GPUs for ordinary differential equations can result in a significant speedup (see \cite{murray2012}). However, a class of problems has been identified, so called stiff problems, where standard integration routines (such as the above mentioned Runge--Kutta methods) are inefficient (see e.g. \cite{Hochbruck2010}).

Exponential integrators are one class of methods that avoid the difficulties of Runge--Kutta method if applied to stiff problems. For such schemes analytical functions (e.g. the exponential function) of large matrices have to be computed. 
%It has long been believed that the efficient computation of such matrix functions present insurmountable difficulties in implementing exponential integrators efficiently. However, due to increased computing power available as well as algorithmic advances this view has changed. 
Exponential integrators and some of their applications are discussed in detail in \cite{Hochbruck2010}. In this paper we will consider a polynomial interpolation scheme to compute the matrix functions; this essentially reduces the problem of efficiently implementing exponential integrators to sparse matrix-vector multiplication as well as computing the nonlinearity of a given differential equation. The computation of matrix-vector multiplications, e.g. by using stencil methods, is usually the most time intensive part of any exponential integrator; thus, an efficient implementation of stencil methods is vital.

\subsection{Research problems \& Results}

In the literature, see section \ref{sec:stencil_methods_stateoftheart}, stencil methods are considered for trivial boundary conditions in the context of a differential operator with constant coefficients (i.e. the Laplacian). Such simplifying assumptions, however, are usually not satisfied in a given application. It is not clear from the literature how much stencil methods can be extended beyond the situation described above while still maintaining an efficient implementation. We propose a method based on the integration of boundary conditions and position-dependent coefficients directly into the CUDA kernel and show that such methods can be applied widely without a significant impact on performance. 

In addition, it has been shown in \cite{Micikevicius2009} that stencil methods can be efficiently parallelized to at least $4$ GPUs. Our objective is to show that such results can be generalized to implementations of exponential integrators for a large class of nonlinearities. 

%\subsection{Organization \& Results}
The remainder of this paper is structured as follows. In section \ref{sec:GPU_arch} a introduction explaining the GPU architecture and the corresponding programming model is given. In addition, we discuss previous work which considers the implementation of stencil methods on GPUs and elaborate on the necessary steps to efficiently implement an exponential integrator (sections \ref{sec:stencil_methods_stateoftheart} and \ref{sec:expint}, respectively).
In section \ref{sec:result} we present our results as summarized below.
\begin{itemize}
	\item Stencil methods that include simple, but non-trivial, boundary conditions, such as those required in many applications, can still be efficiently implemented on GPUs (section \ref{sec:stencil_boundary}). For homogeneous boundary conditions on the C2075  $33.5$ Gflops/s are observed.
	\item Position dependent coefficients (such as a position dependent diffusion) can efficiently be implemented on the GPU provided that the coefficients are not extremely expensive to compute (section \ref{sec:stencil_position}). For a real world example $16$ Gflops/s are observed.
	\item  A wide class of nonlinearities can be computed efficiently on the GPU (section \ref{sec:nonlinearity}).
	\item The parallelization of exponential integrators to multiple GPUs can be conducted very efficiently for discretized differential operators (perfect scaling to at least 4 GPUs) and is mainly limited by the throughput of the PCIe express bus for unstructured matrices (section \ref{sec:multiple-gpus}).
\end{itemize}
%Finally, in section \ref{sec:conclusion}, we discuss the results obtained and provide an outlook for future research. 
Finally, we conclude in section \ref{sec:conclusion}.

%As the computation of matrix-vector multiplication, e.g. by using stencil methods, is usually the most time intensive parts of any exponential integrators a efficient implementation is vital.

%
% put something here that describes the remaining sections
%

\section{BACKGROUND \& MOTIVATION} \label{sec:Background_Motivation}

\subsection{GPU architecture} \label{sec:GPU_arch}

A graphic processing unit (GPU) is a massively parallel computing architecture.
At the time of writing two frameworks to program such systems, namely OpenCL and NVIDIA's CUDA, are in widespread use. In this section we will discuss the hardware architecture as well as the programming model of the GPU architecture using NVIDIA's CUDA (all our implementations are CUDA based). Note, however, that the principles introduced here can, with some change in terminology, just as well be applied to the OpenCL programming model. For a more detailed treatment we refer the reader to \cite{CUDACPG}.

The hardware consists of so called SM (streaming multiprocessors) that are divided into cores.
Each core is (as the analogy suggests) an independent execution unit that shares certain resources (for example shared memory) with other cores, which reside on the same streaming multiprocessor. For example, in case of the C2075, the hardware consists in total of $448$ cores that are distributed across $14$ streaming multiprocessors of $32$ cores each. 

For scheduling, however, the hardware uses the concept of a warp.
A warp is a group of 32 threads that are scheduled to run on the same streaming multiprocessor (but possibly on different cores of that SM). On devices of compute capability 2.0 and higher (e.g. the C2075) each SM consists of 32 cores (matching each thread in a warp to a single core). However, for the C1060, where $240$ cores are distributed across $30$ SM of $8$ core each, even threads in the same warp that take exactly the same execution path are not necessarily scheduled to run in parallel (on the instruction level).

To run a number of threads on a single SM has the advantage that certain resources are shared among those threads; the most notable being the so called shared memory. Shared memory essentially acts as an L1 cache (performance wise) but can be fully controlled by the programmer. Therefore, it is often employed to avoid redundant global memory access as well as to share certain intermediate computations between cores. In addition, devices of compute capability 2.0 and higher are equipped with a non-programmable cache.

The global memory is a RAM (random access memory) that is shared by all SM on the entire GPU. For the C1060 and C2075 GPUs used in this paper the size of the global memory is 4 GB and 6 GB respectively (with memory bandwidth of 102.4 GB/s and 141.7 GB/s respectively), whereas the shared memory is a mere 16~KB for the C1060 and about 50~KB for the C2075. However, this memory is available per SM.

From the programmer some of these details are hidden by the CUDA programming model (most notably the concept of SM, cores, and warps). If a single program is executed on the GPU we refer to this as a grid. The programmer is responsible for subdividing this grid into a number of blocks, whereas each block is further subdivided into threads. A thread in a single block is executed on the same SM and therefore has access to the same shared memory and cache.

GPUs are therefore ideally suited to problems which are compute bound. However, also memory bound problems, such as sparse matrix-vector multiplication, can significantly benefit from GPUs. We will elaborate on this statement in the next section.

\subsection{Stencil methods and matrix-vector products on GPUs} \label{sec:stencil_methods_stateoftheart}

The parallelization of sparse matrix-vector products to GPUs has been studied in some detail. Much research effort in improving the performance of sparse matrix-vector multiplication on GPUs has focused on developing more efficient data structures (see e.g. \cite{Bell2009}  or \cite{Sparse2009}). This is especially important on GPUs as coalesced memory access is of paramount importance if optimal performance is to be achieved. Data structures, such as ELLRT, facilitate coalesced memory access but require additional memory. This is somewhat problematic as on a GPU system memory is limited to a greater extend than on traditional clusters. To remedy this situation a more memory efficient data structure has been proposed, for example, in \cite{A.DziekonskiA.Lamecki2011}. Nevertheless, all such methods are extremely memory bound.

On the other hand, the parallelization of finite difference computations (called stencil methods in this context) to GPUs has been studied, for example, in \cite{Micikevicius2009} and \cite{Berkeley2009}. Even though such methods do not have to store the matrix in memory they are still memory bound; for example, in \cite{Berkeley2009} the flops per byte ratio is computed to be $0.5$ for a seven-point stencil (for double precision computations) which is still far from the theoretical rate of $3.5$ that a C0275 can achieve. In \cite{Berkeley2009} a performance of $36.5$ Gflops/s has been demonstrated for a seven-point stencil on a GTX280.

Both papers mentioned above do not consider boundary conditions in any detail. However, in applications of science and engineering where exponential integrators are applied at least simple boundary conditions have to be imposed (see e.g. \cite{Hochbruck2010}). In addition, in the literature stated above only the discretization of the Laplacian is considered. However, often position-dependent coefficients have to be employed (to model a position-dependent diffusion as in \cite{vazquez2008diffusion}, for example). In this case it is not clear if stencil methods retain their superior performance characteristics (as compared to schemes that store the matrix in memory). We will show in sections \ref{sec:stencil_boundary} and \ref{sec:stencil_position} that for many applications both of these difficulties can be overcome and stencil methods on GPUs can be implemented efficiently.

%To conclude this discussion let us note that more recently, in \cite{Holewinski2012}, an algorithm has been introduced that generates CUDA kernels from a high level description of the finite difference scheme. The performance of the generated code, however, is still a factor of at least $1.5$ away from hand tuned implementations (see \cite{Berkeley2009}).

\subsection{Exponential integrators} \label{sec:expint}
	The step size for the time integration of stiff ordinary differential equations (or the semidiscretization of partial differential equations) is usually limited by a stability condition. In order to overcome this difficulty, implicit schemes are employed that are usually stable for much larger step sizes; however, such schemes have to solve a nonlinear system of equations in each time step and are thus costly in terms of performance.
		 	In many instances the stiffness of the differential equation is located in the linear part only. In this instance, we can write our differential equation as a semilinear problem
		 	\begin{equation}
		 		\frac{\mathrm{d}}{\mathrm{d}t} u(t) + Au(t) = g(u(t)), \label{eq:semilinear}
		 	\end{equation}
		 	where in many applications $A$ is a matrix with large negative eigenvalues and $g$ is a nonlinear function of $u(t)$; it is further assumed that appropriate initial conditions are given. The boundary conditions are incorporated into the matrix $A$.  Since the linear part can be solved exactly, a first-order method, the exponential Euler method, is given by
		 	\begin{equation}
		 		u_{n+1}=\text{e}^{-hA}u_{n}+h\varphi_{1}\left(-hA\right)g(u_{n}), \label{eq:exponential_euler}
		 	\end{equation}
		 	where $\varphi_1$ is an entire function. In \cite{Hochbruck2010} a review of such methods, called exponential integrators,  is given and various methods of higher order are discussed. The main advantage, compared to Runge--Kutta methods, is that an explicit method is given for which the step size is only limited by the nonlinearity.
		 	It has long been believed that the computation of the matrix functions in \eqref{eq:exponential_euler} can not be carried out efficiently. However, if a bound of the field of values of $A$ is known a priori, for example, polynomial interpolation is a viable option. In this case the application of Horner's scheme reduces the problem to repeated matrix-vector products of the form
	 	\begin{equation}
		 		 (\alpha A + \beta I) x, \label{eq:matrix-vector-multiplication}
		 \end{equation}
		 where $A\in\mathbb{K}^{n\times n}$ is a sparse matrix, $I$ is the identity matrix, $x\in\mathbb{K}^n$, and $\alpha, \beta \in \mathbb{K}$ with $\mathbb{K} \in \left\{ \mathbb{R}, \mathbb{C} \right\}$.
		That such a product can be parallelized to small clusters has been shown in \cite{Caliari2004} (for an advection-diffusion equation that is discretized in space by finite differences).
	
		Finally, let us discuss the evaluation of the nonlinearity. In many instances the nonlinearity can be computed pointwise. In this case its evaluation is expected to be easily parallelizable to GPUs. In section \ref{sec:nonlinearity} this behavior is confirmed by numerical experiments.
		If the nonlinearity does include differential operators, such as in Burgers' equation, the evaluation is essentially reduced to sparse-matrix vector multiplication, which we will discuss in some detail in this paper (in the context of stencil methods).
		
\section{RESULTS} \label{sec:result}

	%
 	%	Stencil methods and boundary conditions
 	%
 	\subsection{Stencil methods with boundary conditions} \label{sec:stencil_boundary}

	Let us focus our attention first on the standard seven-point stencil resulting from a discretization of the Laplacian in three dimensions, i.e.
	\begin{align*}
 			(\Delta x)^{2} \left(Au\right)_{i_{x},i_{y},i_{z}}=& -6u_{i_{x},i_{y},i_{z}}\\
 			& +u_{i_{x}+1,i_{y},i_{z}}+u_{i_{x}-1,i_{y},i_{z}}\\
 			& +u_{i_{x},i_{y}+1,i_{z}}+u_{i_{x},i_{y}-1,i_{z}}\\
 			& +u_{i_{x},i_{y},i_{z}+1}+u_{i_{x},i_{y},i_{z}-1},
 	\end{align*} 	
 	where $\Delta x$ is the spacing of the grid points.
	The corresponding matrix-vector product given in \eqref{eq:matrix-vector-multiplication} can then be computed without storing the matrix in memory. For each grid point we have to perform at least $2$ memory operations (a single read and a single store) as well as $10$ floating point operations (6 additions and 2 multiplication for the matrix-vector product as well as single addition and multiplication for the second part of \eqref{eq:matrix-vector-multiplication}).
	
	One could implement a stencil method that employs $8$ memory transactions for every grid point. Following \cite{Berkeley2009} we call this the \textit{naive} method. On the other hand we can try to minimize memory access by storing values in shared memory or the cache (note that the C1060 does not feature a cache but the C2075 does). Since no significant 3D slice fits into the relatively limited  shared memory/cache of both the C1060 and C2075, we take only a 2D slice and iterate over the remaining index. Similar methods have been implemented, for example, in \cite{Berkeley2009} and \cite{Micikevicius2009}. We will call this the \textit{optimized} method. 
	
	To implement boundary conditions we have two options. First, a stencil method can be implemented that considers only grid points that lie strictly in the interior of the domain. Second, we can implement the boundary conditions directly into the CUDA kernel. The approach has the advantage that all computations can be done in a single kernel launch. However, conditional statements have to be inserted into the kernel. Since the kernel is memory bound, we do not expect a significant performance decrease at least for boundary conditions that do not involve the evaluation of complicated functions.
	
	The results of our numerical experiments (for both the naive and optimized method) are given in Table \ref{tab:result_laplacian}. Before we discuss the results let us note that on a dual socket Intel Xeon E5355 system the aggressively hand optimized stencil method implemented in \cite{Berkeley2009} gives $2.5$ Gflops/s.

 \begin{table*}[t]
 		\caption{Timing of a single stencil based matrix-vector computation for a number of implementations and boundary conditions. The corresponding Gflops/s are shown in parentheses. All computations are performed with $n=256^3$. \label{tab:result_laplacian}} 		
		\renewcommand{\arraystretch}{1.3}
		\begin{center}
			\begin{tabular}{|c|c|c|c|c|}		
				\hline
				Device & Boundary & Method & Double& Single \\
				\hline
				\multirow{6}{*}{C1060}
				& \multirow{2}{*}{None} & Stencil (naive) & $13.8$ ms ($12$ Gflops/s) & $8.2$ ms ($20.5$ Gflops/s) \\
				&	& Stencil (optimized) & $7.6$ ms ($22$ Gflops/s) & $8.0$ ms ($21$ Gflops/s) \\
				\cline{2-5}
				& \multirow{2}{*}{\pbox{5cm}{Homogeneous\\  Dirichlet}} &
				Stencil (naive) & $13.4$ ms ($12.5$ Gflops/s) & $7.6$ ms ($22$ Gflops/s) \\
				&	& Stencil (optimized) & $8.8$ ms ($19$ Gflops/s) & $9.2$ ms ($18$ Gflops/s) \\
				\cline{2-5}
				& \multirow{1}{*}{$z(1-z)xy$}
				% & Stencil (naive) & a ms & $15$ ms \\
				& Stencil (optimized) & $36$ ms ($4.5$ Gflops/s) & $39$ ms ($4.5$ Gflops/s) \\				
				\cline{2-5}
				& \multirow{1}{*}{$\sin(\pi z)\exp(-x y)$}
				% & Stencil (naive) & a ms & b \\
				& Stencil (optimized) & $54$ ms ($3$ Gflops/s) & $56$ ms ($3$ Gflops/s) \\										
				\hline

				\multirow{8}{*}{C2075}
				& \multirow{2}{*}{None} &
				Stencil (naive) & $5.5$ ms ($30.5$ Gflops/s) & $3.1$ ms ($54$ Gflops/s) \\
				&	& Stencil (optimized) & $4.3$ ms ($39$ Gflops/s) & $2.9$ ms ($58$ Gflops/s) \\
				\cline{2-5}
				& \multirow{2}{*}{\pbox{5cm}{Homogeneous\\  Dirichlet}} &
				Stencil (naive) & $6$ ms ($28$ Gflops/s) & $3.5$ ms ($48$ Gflops/s) \\
				&	& Stencil (optimized) & $5$ ms ($33.5$ Gflops/s) & $3.9$ ms ($43$ Gflops/s) \\
				\cline{2-5}
				& \multirow{2}{*}{$z(1-z)xy$} &
				Stencil (naive) & $12.3$ ms ($13.5$ Gflops/s) & $6.9$ ms ($24$ Gflops/s) \\
				&	& Stencil (optimized) & $7$ ms ($24$ Gflops/s) & $6.0$ ms ($28$ Gflops/s) \\				
				\cline{2-5}				
				& \multirow{2}{*}{$\sin(\pi z)\exp(-x y)$} &
				Stencil (naive) & $14.3$ ms ($11.5$ Gflops/s) & $13.8$ ms ($12$ Gflops/s) \\
				&	& Stencil (optimized) & $9.7$ ms ($17.5$ Gflops/s)  & $6.8$ ms ($24.5$ Gflops/s) \\							
				\hline
			\end{tabular}	
			
		\par\end{center} 	
 	\end{table*}
		%; that is only slightly better than the $1.7$ Gflops/s we get for a CSR based implementation on a dual socket Intel Xenon E5620.
	
	In \cite{Berkeley2009} a double precision performance of $36.5$ \mbox{Gflops/s} is reported for a GTX280 of compute capability 1.3. However, the theoretical memory bandwidth of the GTX280 is $141.7$ GB/s and thus, as we have a memory bound problem, it has to be compared mainly to the C2075 (which has the same memory bandwidth as the GTX280). Note that the C1060 (compute capability 1.1) has only a memory bandwidth  of $102.4$ GB/s. In our case we  get $39$ Gflops/s for no boundary conditions and $33.5$ Gflops/s for homogeneous Dirichlet boundary conditions. Since we do not solve exactly the same problem, a direct comparison is difficult. However, it is clear that the implemented method is competitive especially since we do not employ any tuning of the kernel parameters.
	
	We found it interesting that for the C2075 (compute capability 2.0) there is only a maximum of $30$\% performance decrease if the naive method is used instead of the optimized method for none or homogeneous boundary conditions (both in the single and double precision case). Thus, the cache implemented on a C2075 works quite efficiently in this case. However, we can get a significant increase in performance for more complicated boundary conditions by using the optimized method. In the single precision case the expected gain is offset, in some instances, by the additional operations that have to be performed in the kernel (see Table \ref{tab:result_laplacian}).
	
	Finally, we observe that the performance for homogeneous Dirichlet boundary conditions is at most $10$\% worse than the same computation which includes no boundary conditions at all. This difference completely vanishes if one considers the optimized implementation. This is no longer true if more complicated boundary conditions are prescribed. For example, if we set
	\[
		f(x,y,z) = z(1-z)xy,
	\]
	for $(x,y,z) \in \partial([0,1]^3)$ or
	\[
		f(x,y,z) = \sin(\pi z)\exp(-x y),
	\]
	for $(x,y,z) \in \partial([0,1]^3)$,
	the performance is decreased by a factor of about $2$ for the C2075 and by a factor of $5$-$7$ for the C1060. Thus, in this case it is clearly warranted to perform the computation of the boundary conditions in a separate kernel launch. Note, however, that the direct implementation is still faster by a factor of~$3$ as compared to CUSPARSE and about $40$ \% better than the ELL format (see \cite{Bell2008}). The memory requirements are an even bigger factor in favor of stencil methods; a grid of dimension $512^3$ would already require $10$ GB in the storage friendly CSR format. Furthermore, the implementation of such a kernel is straight forward and requires no division of the domain into the interior and the boundary.

 \subsection{Stencil methods with a position-dependent coefficient} \label{sec:stencil_position}

	Let us now discuss the addition of a position-dependent diffusion coefficient, i.e. we implement the discretization of $D(x,y,z) \Delta u$ as a stencil method (this is the diffusive part of $\nabla \cdot (D \nabla u )$ ). Compared to the previous section we expect that the direct implementation of the position-dependent diffusion coefficient in the CUDA kernel, for a sufficiently complicated $D$, results in an compute bound problem. For the particular choice of $D(x,y,z)=1/\sqrt{1+x^2+y^2}$, taken from \cite{vazquez2008diffusion}, the results are shown in Table \ref{tab:result_laplacian_nonconstdiffusion}.

 		 \begin{table}[h]
 		\caption{Timing of a single stencil based matrix-vector computation for a position dependent diffusion coefficient given by $D(x,y,z)=1/\sqrt{1+x^2+y^2}$. All computations are performed with  $n=256^3$. \label{tab:result_laplacian_nonconstdiffusion}}
 		\renewcommand{\arraystretch}{1.3}
 				
		\begin{center}
			
			\begin{tabular}{|c|c|c|}
			\multicolumn{3}{c}{\textbf{Double precision} } \\
				\hline
				Device & Method & Time \\
				\hline
				\multirow{2}{*}{C1060}
				& Stencil (naive) & $37$  ms ($4.5$ Gflops/s) \\
				& Stencil (optimized) & $42$ ms ($4$ Gflops/s) \\								
				\hline

				\multirow{2}{*}{C2075}
				& Stencil (naive) & $10.7$ ms ($15.5$ Gflops/s)  \\
				& Stencil (optimized) & $10.5$ ms ($16$ Gflops/s) \\									
				\hline
			\end{tabular}	
			
			\medskip
				
			\begin{tabular}{|c|c|c|}
			\multicolumn{3}{c}{\textbf{Single precision} } \\
				\hline
				Device & Method & Time \\
				\hline
				\multirow{2}{*}{C1060}
				& Stencil (naive) &$37$ ms ($4.5$ Gflops/s) \\
				& Stencil (optimized) & $45$ ms ($3.5$ Gflops/s) \\								
				\hline

				\multirow{2}{*}{C2075}
				& Stencil (naive) & $10.2$ ms ($16.5$ Gflops/s) \\
				& Stencil (optimized) &  $10.3$ ms ($16$ Gflops/s) \\									
				\hline
			\end{tabular}				
			
		\par\end{center} 	
 	\end{table}	
	
	Thus, a performance of $16$ Gflops/s can be achieved for this particular position-dependent diffusion coefficient. This is a significant increase in performance as compared to a matrix-based implementation. In addition, the same concerns regarding storage requirements, as raised above, still apply equally to this problem. No significant difference between the naive and optimized implementation can be observed; this is due to the fact that this problem is now to a large extend compute bound.
	
	Finally, let us note that the results obtained in Tables \ref{tab:result_laplacian} and \ref{tab:result_laplacian_nonconstdiffusion} are (almost) identical for the $n=512^3$ case. Thus, for the sake of brevity, we choose to omit those results.

 	%
 	%	Evaluating the nonlinearity on a GPU
 	%
 	\subsection{Evaluating the nonlinearity on a GPU \label{sec:nonlinearity} }

	For an exponential integrator, usually the most time consuming part is evaluating the exponential and $\varphi_1$ function. Fortunately, if the field of values of $A$ can be estimated a priori, we can employ polynomial interpolation to reduce that problem to matrix-vector multiplication; a viable possibility is interpolation at Leja points (see \cite{Hochbruck2010}). Then, our problem reduces to the evaluation of a series of matrix-vector products of the form given in \eqref{eq:matrix-vector-multiplication} and discussed in the previous section and the evaluation of the nonlinearity for a number of intermediate approximations. In this section we will be concerned with the efficient evaluation of the nonlinearity on a GPU.
	
	Since the nonlinearity is highly problem dependent, let us -- for the sake of concreteness -- take a simple model problem, namely the reaction-diffusion equation modeling combustion in three dimensions (see \cite[p. 439]{hundsdorfer2007})
	\begin{equation}  \label{eq:combustion}
		u_{t}=\Delta v + g(u)
	\end{equation}
%	\begin{equation} \label{eq:combustion}
% 			\begin{cases}
%u_{t}=\Delta v + g(u)\\
%u\vert_{\{t=0\}}=1\\
%u\vert_{\partial \Omega \backslash \{(x,y,z)=1\}}=0, \\
%\partial_x u\vert_{\{(x,y,z)=1\}}=1,
%\end{cases}
% 	\end{equation}
 		with nonlinearity
 	\begin{equation*}
 		g(u)= \frac{1}{4} (2-u) \mathrm{e}^{20 \left(1-\frac{1}{u} \right)}
 	\end{equation*}
 	and appropriate boundary conditions as well as an initial condition.
% 		It is well known that the solution takes values in $[1,2]$ (see \cite[p. 439]{hundsdorfer2007}) and thus the nonlinearity can be large but is bounded.

 		In addition to the discretization of the Laplacian which can be conducted by stencil methods (as described in section \ref{sec:stencil_boundary}) the parallelization of the nonlinearity can be conducted pointwise on the GPU. That is (in a linear indexing scheme) we have to compute
 		\begin{equation} \label{eq:discrete-nonlinearity}
 			\frac{1}{4} (2-u_i) \mathrm{e}^{20 \left(1-\frac{1}{u_i} \right)}, \qquad 0\leq i < n.
 		\end{equation}
		This computation requires only two memory operations per grid point (one read and one store); however, we have to perform a single division and a single exponentiation. Since those operations are expensive,
		% we expect to perform at least the equivalent of $100$ arithmetic instructions. Therefore, the  flops per byte ratio is about $6$ (for double precision) and thus%
		 the problem is expected to be compute bound. The results of our numerical experiments are shown in Table \ref{result_nonlinearity}.
		
		 \begin{table}[h]
 		\caption{Timing of a single computation of the nonlinearity given in \eqref{eq:discrete-nonlinearity}. Results for both full precision computations as well as the fast math routines implemented in the GPU are listed. As a reference a comparison to a dual socket Intel Xenon E5620 setup is provided.  \label{result_nonlinearity}}
	\renewcommand{\arraystretch}{1.3} 				
 				
		\begin{center}
			\begin{tabular}{|c|c|c|c|}
				\multicolumn{4}{c}{\textbf{Double precision} } \\
				\hline
				Device & Method & $n=256^3$ & $n=512^3$ \\
				\hline
				2x Xenon E5620 & OpenMP & $480$ ms & $4$ s \\
				\hline
				
				\multirow{2}{*}{C1060} & Full precision & $14.6$ ms & $120$ ms  \\
					\cline{2-4}
					& Fast math & $6.9$ ms & $55$ ms \\
				\hline
				
				\multirow{2}{*}{C2075} & Full precision & $4.2$ ms & $33$ ms \\
					\cline{2-4}
					& Fast math & $2.4$ ms & $20$ ms \\
				\hline
			\end{tabular}		

			\medskip
		
			\begin{tabular}{|c|c|c|c|}
				\multicolumn{4}{c}{\textbf{Single precision} } \\
				\hline
				Device & Method & $n=256^3$ & $n=512^3$ \\
				\hline
				2x Xenon E5620 & OpenMP & $515$ ms & $4$ s \\
				\hline
				
				\multirow{2}{*}{C1060} & Full precision & $15.4$ ms & $120$ ms \\
					\cline{2-4}
					& Fast math & $7.6$ ms & $61$ ms \\
				\hline
				
				\multirow{2}{*}{C2075} & Full precision & $2.6$ ms & $34$ ms \\
					\cline{2-4}
					& Fast math & $1.6$ ms & $19$ ms \\
				\hline
			\end{tabular}
		\par\end{center}		

 	\end{table}
	As expected, the GPU has a significant advantage over our CPU based system in this case. Fast math routines can be employed if precision is not critical and the evaluation of the nonlinearity contributes significantly to the runtime of the program. Let us duly note that the speedups observed here can not be extended to the entire exponential integrator as the sparse-matrix vector multiplication is usually the limiting factor.
		
		The nonlinearity of certain semi-linear PDEs resemble more the performance characteristics of the stencil methods discussed in sections \ref{sec:stencil_boundary} and \ref{sec:stencil_position}. For example, Burgers' equation, where $g(\boldsymbol{u})=(\boldsymbol{u}\cdot\nabla) \boldsymbol{u}$, falls into this category. Such nonlinearities can be efficiently implemented by the methods discussed in sections \ref{sec:stencil_boundary} and \ref{sec:stencil_position}.
		
		If we combine sections \ref{sec:stencil_boundary}, \ref{sec:stencil_position} and \ref{sec:nonlinearity} we have all ingredients necessary to conduct an efficient implementation of exponential integrators on a single GPU. The specific performance characteristics depend on the form of the linear as well as the nonlinear part of the differential equation under consideration. In the next section we will turn our attention to the parallelization of exponential integrators to multiple GPUs.

 	%
 	%	Exponential integrators on multiple GPU's
 	%
\subsection{Multiple GPU implementation of exponential integrators \label{sec:multiple-gpus} }

	If we consider the problem introduced in \eqref{eq:combustion} to be solved with an exponential integrator, we have at least two possibilities to distribute the workload to multiple GPUs. First, one could compute the different matrix functions on different GPUs. However, since even for higher order schemes we only have to evaluate a small number of distinct matrix functions, this approach is not very flexible and depends on the method under consideration. However, if we are able to implement a parallelization of the matrix-vector product and the nonlinearity onto multiple GPUs, a much more flexible approach would result.
	
	Such an undertaking however is limited by the fact that in the worst case we have to transfer
	\begin{equation} \label{eq:worst-case}
		(m-1)n
	\end{equation}
	floating point numbers over the relatively slow PCIe bus ($m$ is the number of GPUs whereas $n$ is, as before, the problem size). However, in the case of differential operators only a halo region has to be updated after every iteration and thus the actual memory transfer is a tiny fraction of the value given by \eqref{eq:worst-case}. Such a procedure  was suggested in \cite{Caliari2004} for use on a cluster, where parallelization is mainly limited by the interconnection between different nodes.
	For performance reasons on a GPU it is advantageous to first flatten the halo regions in memory and copy it via a single call to {\tt cudaMemcpy} to the device. Then the vector is updated by using that information in a fully parallelized way on the GPU. As can be seen from the results given in Table \ref{tab:compustion-performance}, the problem in \eqref{eq:combustion} shows good scaling behavior (at least) up to 4 GPUs.
	\begin{table}[H]
 		\caption{Performance comparison for the combustion model discussed in section \ref{sec:nonlinearity} for a single time step using $40$ matrix-vector products (a tolerance of $\text{tol}=10^{-4}$ was prescribed for a time step of size $10^{-4}$). A finite difference discretization with $n=256^3$ has been used. \label{tab:compustion-performance} }
	\renewcommand{\arraystretch}{1.3} 		 		
 		
		\begin{center}
			\begin{tabular}{|c|c|c|c|}
				\multicolumn{4}{c}{\textbf{Double precision} } \\
				\hline
				Device & Method &  \pbox{3cm}{\vspace{0.2cm}Number \\units} & Time \\[0.5cm]
				\hline
				2x Xenon E5620 &  CSR/OpenMP & 2 & $9.5$ s   \\
				\hline
				\multirow{2}{*}{C1060} & \multirow{2}{*}{\pbox{3cm}{Stencil\\hom. Dirichlet}} & 1 &  $1.5$ s    \\
					\cline{3-4}
					%& & 2 & a  &$730$ ms  \\
					\cline{3-4}
					& & 4 & $320$ ms   \\	
				\hline
				C2075 & \pbox{3cm}{\vspace{0.1cm}Stencil\\hom. Dirichlet} & 1 & $1.2$ s  \\[0.2cm]
				\hline
			\end{tabular}	
		
			\medskip
		
			\begin{tabular}{|c|c|c|c|}
				\multicolumn{4}{c}{\textbf{Single precision} } \\
				\hline
				Device & Method &  \pbox{3cm}{\vspace{0.2cm}Number \\units} & Time \\[0.5cm]
				\hline
				2x Xenon E5620 &  CSR/OpenMP & 2  & $5.6$ s  \\
				\hline
				\multirow{2}{*}{C1060} & \multirow{2}{*}{\pbox{3cm}{Stencil\\hom. Dirichlet}} & 1 &  $1.2$ s  \\
					\cline{3-4}
					%& & 2 & a  &$730$ ms  \\
					\cline{3-4}
					& & 4   & $540$ ms  \\	
				\hline
				C2075 & \pbox{3cm}{\vspace{0.1cm}Stencil\\hom. Dirichlet} & 1 & $210$ ms \\[0.2cm]
				\hline
			\end{tabular}

		\par\end{center}

		% $630$ ms C1060 1 GPU: $360$ ms
		% $130$ ms,  $130$ ms
	\end{table}	

	Let us now discuss a different example. In certain discrete quantum systems, for example, the solution of (see, e.g., \cite{Raedt08})
	\[
		\partial_t \psi = H(t)\psi
	\]
	is to be determined, where $\psi$ is a vector with complex entries in a high dimensional vector space and $H(t)$ a Hermitian matrix. Such equations are efficiently solved by using Magnus integrators. In this paper we will use the example of a two spin system in a spin bath. In this case $H(t)$ is independent of time and thus we can, in principle, take arbitrarily large time steps. The matrix $H$ is generated beforehand and stored in the generic CSR format; for $21$ spins this yields a vector with $n=2^{21}$ complex entries and a matrix with approximately $83.9\cdot 10^6$ non-zero complex entries (the storage requirement is about $2$ GB in the double precision and $1$ GB in the single precision case). This gives a sparsity of $2\cdot 10^{-5}.$ Note, however, that such quantum systems couple every degree of freedom with every other degree of freedom. Thus, we are in the worst case and have to transfer \mbox{$(m-1)n$} floating point numbers over the PCIe bus after each iteration. 
	
	The results of our numerical experiments are shown in Table \ref{tab:results-spin-bath}. The implementation used is based on the code given in \cite{NVIDIA:2008}. However, we have found that for the problem under consideration using a full warp for every row of our matrix results in a performance reduction. Therefore, we use only four threads per row which results, for the specific problem under consideration, in a performance increase of approximately $50$\%, as compared against the CUSPARSE library (see \cite{Cusparse}). Apart from this consideration the code has been adapted to compute the  problem stated in \eqref{eq:matrix-vector-multiplication}, which includes an additional term as compared to the sparse matrix-vector multiplication considered in \cite{NVIDIA:2008}.
	 \begin{table}[h]
 		\caption{Performance comparison for a system with $21$ spins. Integration is performed up to $t=10$ with a tolerance of $\text{tol}=10^{-5}$. \label{tab:results-spin-bath}}
 			\renewcommand{\arraystretch}{1.3} 	
		\begin{center}
		
			\begin{tabular}{|c|c|c|c|}
				\multicolumn{4}{c}{\textbf{Double precision} } \\
				\hline
				Device & Method & \pbox{3cm}{\vspace{0.2cm}Number \\units} & Time \\[0.5cm]
				\hline
				2x Xenon E5620 &  CSR/OpenMP & 2 & $46$ s \\
				\hline
				\multirow{3}{*}{C1060} & \multirow{3}{*}{CSR} & 1 & $23$ s \\
					\cline{3-4}
					& & 2 & $15$ s \\
					\cline{3-4}
					& & 4 & $15$ s \\	
				\hline
				C2075 & CSR & 1 & $7.7$ s \\
				\hline
			\end{tabular}		
		
			\medskip		
		
			\begin{tabular}{|c|c|c|c|}
				\multicolumn{4}{c}{\textbf{Single precision} } \\
				\hline
				Device & Method & \pbox{3cm}{\vspace{0.2cm}Number \\units} & Time \\[0.5cm]
				\hline
				2x Xenon E5620 &  CSR/OpenMP & 2 & $44$ s \\
				\hline
				\multirow{3}{*}{C1060} & \multirow{3}{*}{CSR} & 1 & $15$ s \\
					\cline{3-4}
					& & 2  & $10$ s \\
					\cline{3-4}
					& & 4 & $7.5$ s \\	
				\hline
				C2075 & CSR & 1 &  $4$ s \\
				\hline
			\end{tabular}			
			
		\par\end{center}		
	\end{table}
	Clearly the scaling behavior in this case is limited by the overhead of copying between the different GPUs. For two GPUs a speedup of about $1.5$ can be observed. For any additional GPU no performance gain can be observed. In total a speedup of $3$ for double precision and $6$ for single precision as compared to a dual socket Xenon configuration is achieved on four C1060 graphic processing units. This is only about $50$\% better than the speedup of $2$ (double precision) and $3$ (single precision) achieved with a single C1060 GPU. It should be noted that a GPU centric data format (as discussed in section \ref{sec:stencil_methods_stateoftheart}) could be employed instead of the CSR format. However, also in this case the overhead of copying between different GPUs would persist.
	
	Thus, in this instance the speedups that are achievable in both single and multi GPU configurations are a consequence of an unstructured matrix that makes coalesced memory access as well as parallelizability between different GPUs difficult. The dramatically better performance of the C2075 as shown in Table \ref{tab:results-spin-bath} is thus expected.

 	%
 	%	Conclusion
 	% 	
\section{CONCLUSION \& OUTLOOK \label{sec:conclusion} }
	We have shown that exponential integrators can be efficiently implemented on graphic processing units. For many problems, especially those resulting from the spatial discretization of partial differential equations, this is true for both single and multi GPU setups.
	
	In addition, we have considered stencil based implementations that go beyond periodic boundary conditions and constant diffusion coefficients. Such problems can not be handled by implementations based on the fast Fourier transform, for example. Moreover, section \ref{sec:stencil_boundary} shows that for non-homogeneous boundary conditions the code handling the interior as well as the boundary of the domain has to be separated if optimal performance is to be achieved. However, for homogeneous or piecewise constant boundary condition an implementation directly into the CUDA kernel does not result in any significant performance decrease.

The results presented in this paper show that exponential integrators, for many realistic settings, can be efficiently implemented on GPUs with significant speedups compared to more traditional implementations. Therefore, GPU computing provides a viable way to increase the efficiency of simulations in which exponential integrators are employed.
	The implementation of exponential integrators on the current generation of GPUs would conceivably result in a further performance increase of our memory bound stencil implementation, as compared to the C2075, as the Kepler architecture offers a memory throughput of up to $250$ GB/s. Furthermore, such an implementation is expected to be relatively straightforward as the cache implemented on the newer generations of GPUs works quite well in the case of stencil methods.
	%(For the C2075 this has been demonstrated in some detail in this paper).

% bibliography
	
	\bibliographystyle{IEEEtran}
	\bibliography{exponential-integrators,GPU,gpu_2}

% that's all folks
\end{document}